\documentclass{article}%
\usepackage{amsfonts}
\usepackage{graphicx}
\usepackage{amsmath}
\usepackage[francais]{babel}
\usepackage{amssymb}%
\newtheorem{theorem}{Theorem}

\newtheorem{proposition}[theorem]{Proposition}

\newtheorem{theoreme}[theorem]{Théorème}
\newtheorem{corollaire}[theorem]{Corollaire}

\newenvironment{demonstration}[1][D\'emonstration]{\textbf{#1:} }
{\ \rule{0.5em}{0.5em}}
\begin{document}

\title{Structures de Jacobi sur une vari\'{e}t\'{e} des points proches.}
\author{Basile Guy Richard BOSSOTO\\Universit\'{e} Marien NGOUABI, Facult\'{e} des Sciences,\\D\'{e}partement de Math\'{e}matiques\\B.P.69 - Brazzaville, Congo\\E-mail: bossotob@yahoo.fr}
\date{}
\maketitle

\begin{abstract}
On consid\`{e}re une alg\`{e}bre locale $A$ (au sens d'Andr\'{e} Weil), $M$
une vari\'{e}t\'{e} lisse paracompacte et $M^{A}$ la vari\'{e}t\'{e} des
points proches de $M$ d'esp\`{e}ce $A$. Dans ce travail, nous d\'{e}finissons
et \'{e}tudions les notions de $A$-structures de Jacobi sur $M^{A}$.

\textbf{Abstract : }We consider a local algebra $A$ (in the sense of Andr\'{e}
Weil), $M$ a smooth paracompact manifold and $M^{A}$ the manifold of
infinietly near points on $M$ of kind $A$. In this paper, we define and study
the notions of $A$-Jacobi structures on $M^{A}$.

\end{abstract}

\begin{quote}

\end{quote}

\section{Introduction}

On consid\`{e}re une alg\`{e}bre locale $A$ (au sens d'Andr\'{e} Weil),
c'est-\`{a}-dire une alg\`{e}bre r\'{e}elle $A$ commutative unitaire de
dimension finie sur $%
\mathbb{R}
$, ayant un id\'{e}al maximal unique de codimension $1$ sur $%
\mathbb{R}
$, $M$ une vari\'{e}t\'{e} lisse paracompacte et $M^{A}$ la vari\'{e}t\'{e}
des points proches de $M$ d'esp\`{e}ce $A$ \cite{wei}.

L'ensemble, $C^{\infty}(M^{A},A)$, des fonctions sur $M^{A}$ \`{a} valeurs
dans $A$ est une alg\`{e}bre commutative unitaire sur $A$. En notant,
$C^{\infty}(M)$, l'alg\`{e}bre des fonctions num\'{e}riques de classe
$C^{\infty}$ sur $M$, alors pour $f\in C^{\infty}(M)$, l'application
\[
f^{A}:M^{A}\longrightarrow A,\xi\longmapsto\xi(f),
\]
est de classe $C^{\infty}$ et l'application%
\[
C^{\infty}(M)\longrightarrow C^{\infty}(M^{A},A),f\longmapsto f^{A},
\]
est un homomorphisme d'alg\`{e}bres r\'{e}elles.

Il y a \'{e}quivalence entre les assertions suivantes \cite{bos}:

\begin{enumerate}
\item $X:C^{\infty}(M^{A})\longrightarrow C^{\infty}(M^{A})$ est un champ de
vecteurs sur $M^{A}$ ;

\item $X:C^{\infty}(M)\longrightarrow C^{\infty}(M^{A},A)$ est une application
lin\'{e}aire v\'{e}rifiant%

\[
X(fg)=X(f)\cdot g^{A}+f^{A}\cdot X(g)
\]
pour tous $f$ et $g$ dans $C^{\infty}(M)$.
\end{enumerate}

Ainsi l'ensemble, $\mathfrak{X}(M^{A})$, des champs de vecteurs sur $M^{A}$
est un $C^{\infty}(M^{A},A)$-module \cite{bos}.

Lorsque $X$ \ est un champ de vecteurs sur $M^{A}$, consid\'{e}r\'{e} comme
d\'{e}rivation de $C^{\infty}(M)$ dans $C^{\infty}(M^{A},A)$, alors il existe
une d\'{e}rivation et une seule \cite{bos}%
\[
\widetilde{X}:C^{\infty}(M^{A},A)\longrightarrow C^{\infty}(M^{A},A)
\]
telle que:

\begin{enumerate}
\item $\widetilde{X}$ est $A$-lin\'{e}aire;

\item $\widetilde{X}\left[  C^{\infty}(M^{A})\right]  \subset C^{\infty}%
(M^{A})$;

\item $\ \widetilde{X}(f^{A})=X(f)$ pour tout $f\in C^{\infty}(M)$.
\end{enumerate}

Lorsque $R$ est un anneau commutatif unitaire, d'\'{e}l\'{e}ment unit\'{e}
$1_{R}$ , et lorsque $E$ est un $R$-module, une application lin\'{e}aire
\[
\delta:R\longrightarrow E
\]
est un op\'{e}rateur diff\'{e}rentiel d'ordre $\leq1$ si, pour tous $a$ et $b$
dans $R$,%
\[
\delta(ab)=\delta(a)\cdot b+a\cdot\delta(b)-ab\cdot\delta(1_{R})\text{.}%
\]
Lorsque $\delta(1_{R})=0$, on a la notion usuelle de d\'{e}rivation de
$R\ $\ dans $E$.

Ainsi une application lin\'{e}aire
\[
\delta:R\longrightarrow E
\]
est un op\'{e}rateur diff\'{e}rentiel d'ordre $\leq1$ si et seulement si
l'application
\[
R\longrightarrow E,a\mapsto\delta(a)-a\cdot\delta(1_{R}),
\]
est une d\'{e}rivation.

Dans toute la suite, $A$ d\'{e}signe une alg\`{e}bre locale (au sens
d'Andr\'{e} Weil), $M$ une vari\'{e}t\'{e} lisse paracompacte, $M^{A}$ la
vari\'{e}t\'{e} des points proches de $M$ d'esp\`{e}ce $A$, $C^{\infty}(M)$
l'alg\`{e}bre des fonctions num\'{e}riques de classe $C^{\infty}$ sur $M$ et
d'\'{e}l\'{e}ment-unit\'{e} $1_{C^{\infty}(M)}$, $\mathfrak{X}(M)$
l'alg\`{e}bre de Lie r\'{e}elle des champs de vecteurs sur $M$ et
$\mathcal{D}(M)$ l'alg\`{e}bre de Lie r\'{e}elle des op\'{e}rateurs
diff\'{e}rentiels d'ordre $\leq1$ de $C^{\infty}(M)$ dans $C^{\infty}(M)$. Le
terme "op\'{e}rateur diff\'{e}rentiel" signifiera "op\'{e}rateur
diff\'{e}rentiel d'odre $\leq1$".

\section{Structure de $A$-alg\`{e}bre de Lie-Rinehart sur $\mathcal{D}(M^{A}%
)$}

\begin{proposition}
Les assertions suivantes sont \'{e}quivalentes:

\begin{enumerate}
\item $X:C^{\infty}(M^{A})\longrightarrow C^{\infty}(M^{A})$ est un
op\'{e}rateur diff\'{e}rentiel;

\item $X:C^{\infty}(M)\longrightarrow C^{\infty}(M^{A},A)$ est une application
$%
\mathbb{R}
$-lin\'{e}aire v\'{e}rifiant%
\[
X(fg)=X(f)\cdot g^{A}+f^{A}\cdot X(g)-f^{A}\cdot g^{A}\cdot X(1_{C^{\infty
}(M)})
\]

pour tous $f,g\in C^{\infty}(M)$.
\end{enumerate}
\end{proposition}

\begin{demonstration}
Soit $(a_{\alpha})_{\alpha\in I}$ une base de $A$ et soit $(a_{\alpha}^{\ast
})_{\alpha\in I}$ la base duale de la base $(a_{\alpha})_{\alpha\in I}$.

$/\Longrightarrow$ Comme $X$ est un op\'{e}rateur diff\'{e}rentiel, alors
l'application
\[
C^{\infty}(M^{A})\longrightarrow C^{\infty}(M^{A}),\varphi\mapsto
X(\varphi)-\varphi\cdot X(1_{C^{\infty}(M^{A})})
\]
est une d\'{e}rivation, donc un champ de vecteurs sur $M^{A}$. Compte tenu de
\cite{bos}, l'application
\[
Y:C^{\infty}(M)\longrightarrow C^{\infty}(M^{A},A),f\longmapsto\left(
{\textstyle\sum\limits_{\alpha\in I}}
\left[  X(a_{\alpha}^{\ast}\circ f^{A})\right]  \cdot a_{\alpha}\right)
-f^{A}\cdot X(1_{C^{\infty}(M^{A})})
\]
v\'{e}rifie
\[
Y(fg)=Y(f)\cdot g^{A}+f^{A}\cdot Y(g)\text{.}%
\]
On v\'{e}rifie que l'application
\[
C^{\infty}(M)\longrightarrow C^{\infty}(M^{A},A),f\longmapsto Y(f)+f^{A}\cdot
X(1_{C^{\infty}(M^{A})})
\]
r\'{e}pond \`{a} la question.

$\Longleftarrow/$ Soit $X:C^{\infty}(M)\longrightarrow C^{\infty}(M^{A},A)$
une application $%
\mathbb{R}
$-lin\'{e}aire v\'{e}rifiant%
\[
X(fg)=X(f)\cdot g^{A}+f^{A}\cdot X(g)-f^{A}\cdot g^{A}\cdot X(1_{C^{\infty
}(M)})
\]
pour tous $f,g\in C^{\infty}(M)$. L'application
\[
Z:C^{\infty}(M)\longrightarrow C^{\infty}(M^{A},A),f\longmapsto X(f)-f^{A}%
\cdot X(1_{C^{\infty}(M)})
\]
est un champ de vecteurs sur $M^{A}$ consid\'{e}r\'{e} comme d\'{e}rivation de
$C^{\infty}(M)$ dans $C^{\infty}(M^{A},A)$. Compte tenu de \cite{bos}, il
existe une d\'{e}rivation et une seule
\[
\overline{Z}:C^{\infty}(M^{A})\longrightarrow C^{\infty}(M^{A})
\]
telle que
\[
Z(f)=%
{\textstyle\sum\limits_{\alpha\in I}}
\overline{Z}(a_{\alpha}^{\ast}\circ f^{A})\cdot a_{\alpha}%
\]
L'application
\[
C^{\infty}(M^{A})\longrightarrow C^{\infty}(M^{A}),\varphi\longmapsto
\overline{Z}(\varphi)+\varphi\cdot X(1_{C^{\infty}(M)}),
\]
est un op\'{e}rateur diff\'{e}rentiel.
\end{demonstration}

L'ensemble, $\mathcal{D}(M^{A})$, des op\'{e}rateurs diff\'{e}rentiels de
$C^{\infty}(M^{A})$ dans $C^{\infty}(M^{A})$ consid\'{e}r\'{e}s comme
op\'{e}rateurs diff\'{e}rentiels de $C^{\infty}(M)\longrightarrow C^{\infty
}(M^{A},A)$ est un $C^{\infty}(M^{A},A)$-module.

\begin{proposition}
Si $X\in\mathcal{D}(M^{A})$, consid\'{e}r\'{e} comme op\'{e}rateur
diff\'{e}rentiel de $C^{\infty}(M)\longrightarrow C^{\infty}(M^{A},A)$, alors
il existe un op\'{e}rateur diff\'{e}rentiel et un seul%
\[
\widetilde{X}:C^{\infty}(M^{A},A)\longrightarrow C^{\infty}(M^{A},A)
\]
tel que :

\begin{enumerate}
\item $\widetilde{X}$ est $A$-lin\'{e}aire;

\item Pour tout $\varphi\in C^{\infty}(M^{A})$, $\left[  \widetilde{X}%
(\varphi)-\varphi\cdot X(1_{C^{\infty}(M)})\right]  \in C^{\infty}(M^{A})$;

\item $\widetilde{X}(f^{A})=X(f)$ pour tout $f\in C^{\infty}(M)$.
\end{enumerate}
\end{proposition}

\begin{demonstration}
Comme \textbf{\ }$X:C^{\infty}(M)\longrightarrow C^{\infty}(M^{A},A)$ est un
op\'{e}rateur diff\'{e}rentiel, alors l'application%
\[
Y:C^{\infty}(M)\longrightarrow C^{\infty}(M^{A},A),f\longmapsto X(f)-f^{A}%
\cdot X(1_{C^{\infty}(M)}),
\]
est une d\'{e}rivation. Compte tenu de \cite{bos}, il existe une
d\'{e}rivation et une seule
\[
\widetilde{Y}:C^{\infty}(M^{A},A)\longrightarrow C^{\infty}(M^{A},A)
\]
telle que:

$1$- $\widetilde{Y}$ est $A$-lin\'{e}aire;

$2$- $\widetilde{Y}(\sigma)\in C^{\infty}(M^{A})$ pour tout $\sigma\in
C^{\infty}(M^{A})$;

$3$- $\widetilde{Y}(f^{A})=Y(f)$ pour tout $f\in C^{\infty}(M)$.

$\ $L'application
\[
\widetilde{X}:C^{\infty}(M^{A},A)\longrightarrow C^{\infty}(M^{A}%
,A),\varphi\longmapsto\widetilde{Y}(\varphi)+\varphi\cdot X(1_{C^{\infty}%
(M)}),
\]

est telle que :

$1$- $\widetilde{X}$ est $A$-lin\'{e}aire;

$2$- Pour tout $\varphi\in C^{\infty}(M^{A})$, $\left[  \widetilde{X}%
(\varphi)-\varphi\cdot X(1_{C^{\infty}(M)})\right]  \in C^{\infty}(M^{A})$;

$3$- $\widetilde{X}(f^{A})=X(f)$ pour tout $f\in C^{\infty}(M)$.

Ce qui ach\`{e}ve la d\'{e}monstration.
\end{demonstration}

\begin{proposition}
Pour $\varphi\in C^{\infty}(M^{A},A)$ et pour $X\in\mathcal{D}(M^{A})$, on a
\[
\widetilde{\varphi\cdot X}=\varphi\cdot\widetilde{X}\text{.}%
\]

\end{proposition}

\begin{theoreme}
L'application%
\[
\left[  ,\right]  :\mathcal{D}(M^{A})\times\mathcal{D}(M^{A})\longrightarrow
\mathcal{D}(M^{A}),(X,Y)\longmapsto\widetilde{X}\circ Y-\widetilde{Y}\circ X,
\]
est $A$-bilin\'{e}aire altern\'{e}e et d\'{e}finit une structure de
$A$-alg\`{e}bre de Lie sur $\mathcal{D}(M^{A})$.

De plus, pour $\varphi\in C^{\infty}(M^{A},A)$ et pour $X,Y\in\mathcal{D}%
(M^{A})$, on a%

\[
\lbrack\widetilde{X},\widetilde{Y}]\ =\widetilde{[X,Y]}%
\]

et
\[
\left[  X,\varphi\cdot Y\right]  =\left[  \widetilde{X}(\varphi)-\varphi
\cdot\widetilde{X}(1_{C^{\infty}(M^{A},A)})\right]  \cdot Y+\varphi
\cdot\left[  X,Y\right]  \text{.}%
\]

\end{theoreme}

\begin{demonstration}
Pour $X,Y\in\mathcal{D}(M^{A})$, l'application
\[
\left[  ,\right]  :\mathcal{D}(M^{A})\times\mathcal{D}(M^{A})\longrightarrow
\mathcal{D}(M^{A}),(X,Y)\longmapsto\widetilde{X}\circ Y-\widetilde{Y}\circ X,
\]
est manifestement $%
\mathbb{R}
$-bilin\'{e}aire altern\'{e}e. Pour $f$ et $g$ appartenant \`{a} $C^{\infty
}(M)$, on v\'{e}rifie que%

\[
\left[  X,Y\right]  (fg)=\left[  X,Y\right]  (f)\cdot g^{A}+f^{A}\cdot\left[
X,Y\right]  (g)-f^{A}\cdot g^{A}\cdot\left[  X,Y\right]  (1_{C^{\infty}(M)})
\]

Ainsi $\left[  X,Y\right]  \in\mathcal{D}(M^{A})$.

Pour $a\in A$ et pour $f$ $\in$ $C^{\infty}(M)$, on a:
\begin{align*}
\left[  X,a\cdot Y\right]  (f\ )  &  =\widetilde{X}\left[  a\cdot
Y(f\ )\right]  -\widetilde{a\cdot Y}\left[  X(f\ )\right] \\
&  =a\cdot\widetilde{X}\left[  Y(f\ )\right]  -a\cdot\widetilde{Y}\left[
X(f\ )\right] \\
&  =a\cdot\left[  X,Y\right]  (f\ )\text{.}%
\end{align*}

Ainsi
\[
\left[  X,a\cdot Y\right]  =a\cdot\left[  X,Y\right]  \text{.}%
\]

L'application%
\[
\left[  ,\right]  :\mathcal{D}(M^{A})\times\mathcal{D}(M^{A})\longrightarrow
\mathcal{D}(M^{A}),(X,Y)\longmapsto\widetilde{X}\circ Y-\widetilde{Y}\circ X,
\]
est donc $A$-bilin\'{e}aire altern\'{e}e.

Pour $X,Y\in\mathcal{D}(M^{A})$ et pour $\varphi\in$ $C^{\infty}(M^{A},A)$,
l'application
\[
C^{\infty}(M^{A},A)\longrightarrow C^{\infty}(M^{A},A),\varphi\longmapsto
\lbrack\widetilde{X},\widetilde{Y}](\varphi)=\widetilde{X}\left[
\widetilde{Y}(\varphi)\right]  -\widetilde{Y}\left[  \widetilde{X}%
(\varphi)\right]
\]
est $A$-lin\'{e}aire.

Pour $\sigma\in$ $C^{\infty}(M^{A})$, on v\'{e}rifie que%
\[
\lbrack\widetilde{X},\widetilde{Y}](\sigma)-\sigma\cdot\left[  X,Y\right]
(1_{C^{\infty}(M)})
\]
appartient \`{a} $C^{\infty}(M^{A})$ et que
\[
\lbrack\widetilde{X},\widetilde{Y}](f^{A})=\left[  X,Y\right]  (f)
\]
pour tout $f$ $\in$ $C^{\infty}(M)$. On d\'{e}duit que
\[
\lbrack\widetilde{X},\widetilde{Y}]\ =\widetilde{[X,Y]}\text{.}%
\]

Pour $\varphi\in C^{\infty}(M^{A},A)$ et pour $f$ $\in$ $C^{\infty}(M)$, on a%

\begin{align*}
\left[  X,\varphi\cdot Y\right]  (f\ )  &  =\widetilde{X}\left[  \varphi\cdot
Y(f\ )\right]  -\widetilde{\varphi\cdot Y}\left[  X(f\ )\right] \\
&  =\widetilde{X}(\varphi)\cdot Y(f\ )+\varphi\cdot\widetilde{X}\left[
Y(f\ )\right]  -\varphi\cdot Y(f\ )\cdot\widetilde{X}(1_{C^{\infty}(M^{A}%
,A)})-\varphi\cdot\widetilde{Y}\left[  X(f\ )\right] \\
&  =\widetilde{X}(\varphi)\cdot Y(f\ )-\varphi\cdot Y(f\ )\cdot\widetilde
{X}(1_{C^{\infty}(M^{A},A)})+\varphi\cdot\widetilde{X}\left[  Y(f\ )\right]
-\varphi\cdot\widetilde{Y}\left[  X(f\ )\right] \\
&  =\widetilde{X}(\varphi)\cdot Y(f\ )-\varphi\cdot Y(f\ )\cdot\widetilde
{X}(1_{C^{\infty}(M^{A},A)})+\varphi\cdot\left[  X,Y\right]  (f\ )\\
&  =\left[  \widetilde{X}(\varphi)-\varphi\cdot\widetilde{X}(1_{C^{\infty
}(M^{A},A)})\right]  \cdot Y(f\ )+\varphi\cdot\left[  X,Y\right]  (f\ )
\end{align*}

Ainsi
\[
\left[  X,\varphi\cdot Y\right]  =\left[  \widetilde{X}(\varphi)-\varphi
\cdot\widetilde{X}(1_{C^{\infty}(M^{A},A)})\right]  \cdot Y+\varphi
\cdot\left[  X,Y\right]
\]

D'o\`{u} l'assertion.
\end{demonstration}

Lorsque $\mathcal{D}_{A}\left[  C^{\infty}(M^{A},A)\right]  $ d\'{e}signe le
$C^{\infty}(M^{A},A)$-module des op\'{e}rateurs diff\'{e}rentiels de
$C^{\infty}(M^{A},A)$ dans $C^{\infty}(M^{A},A)$ qui sont $A$-lin\'{e}aires,
l' application
\[
\mathbf{\sim}:\mathcal{D}(M^{A})\longrightarrow\mathcal{D}_{A}\left[
C^{\infty}(M^{A},A)\right]  ,X\longmapsto\widetilde{X},
\]
est $C^{\infty}(M^{A},A)$-lin\'{e}aire et est un morphisme de $A$-alg\`{e}bres
de Lie.

En suivant \cite{ok3}, on a:

\begin{corollaire}
Le couple $(\mathcal{D}(M^{A}),\mathbf{\sim})$ est une $A$-alg\`{e}bre de Lie-Rinehart.
\end{corollaire}

\section{$A$-structures de Jacobi sur $M^{A}$}

\bigskip Une $A$-structure de Jacobi sur $M^{A}$ ou une structure de
$A$-alg\`{e}bre de Jacobi sur $C^{\infty}(M^{A},A)$ est la donn\'{e}e d'une
structure de $A$-alg\`{e}bre de Lie sur $C^{\infty}(M^{A},A)$, de crochet
$\left\{  ,\right\}  $, telle l'application
\[
ad(\varphi):C^{\infty}(M^{A},A)\longrightarrow C^{\infty}(M^{A},A),\psi
\mapsto\left\{  \varphi,\psi\right\}  ,
\]
soit un $\ $op\'{e}rateur diff\'{e}rentiel pour tout $\varphi\in C^{\infty
}(M^{A},A)$.

On va, dans ce qui suit, construire trois $A$-structures de Jacobi sur $M^{A}$.

\subsection{$A$-structure de Jacobi sur $M^{A}$ lorsque le couple
$(\mathcal{D}(M^{A}),\mathbf{\sim})$ admet une structure de $A$-alg\`{e}bre de
Lie-Rinehart-Jacobi symplectique}

On rappelle que le couple $(\mathcal{D}(M^{A}),\mathbf{\sim})$ admet une
structure de $A$-alg\`{e}bre de Lie-Rinehart-Jacobi symplectique s'il existe
une $2$-forme altern\'{e}e et nond\'{e}g\'{e}n\'{e}r\'{e}e
\[
\Omega:\mathcal{D}(M^{A})\times\mathcal{D}(M^{A})\longrightarrow C^{\infty
}(M^{A},A)
\]
telle que $\widetilde{d}\Omega=0$, o\`{u} $\widetilde{d}$ est la
diff\'{e}rentielle de degr\'{e} $+1$ associ\'{e}e \`{a} la repr\'{e}sentation
$\mathbf{\sim}$ \cite{ok3}.

Pour $\varphi\in C^{\infty}(M^{A},A)$, on note $X_{\varphi}$ l'unique
\'{e}l\'{e}ment de $\mathcal{D}(M^{A})$ tel que
\[
i_{X_{\varphi}}\Omega=\widetilde{d}\varphi\text{.}%
\]

Pour tout $\varphi\in C^{\infty}(M^{A},A)$, on v\'{e}rifie que $\theta
_{X_{\varphi}}\Omega=0$ o\`{u} \bigskip\ $\theta_{X_{\varphi}}$ d\'{e}signe la
d\'{e}riv\'{e}e de Lie par rapport \`{a} l'op\'{e}rateur diff\'{e}rentiel
$X_{\varphi}$.

Pour $\varphi,\psi\in C^{\infty}(M^{A},A)$, on pose
\[
\left\{  \varphi,\psi\right\}  =-\Omega(X_{\varphi},X_{\psi})\text{.}%
\]

\begin{proposition}
Si $\varphi,\psi\in C^{\infty}(M^{A},A)$, alors
\[
\left\{  \varphi,\psi\right\}  =\widetilde{X_{\varphi}}(\psi)
\]
et
\[
\left[  X_{\varphi},X_{\psi}\right]  =X_{\left\{  \varphi,\psi\right\}
}\text{.}%
\]

\end{proposition}

\begin{demonstration}
\bigskip\ \bigskip Pour $\varphi,\psi\in C^{\infty}(M^{A},A)$, on a
\begin{align*}
\{\varphi,\psi\}  &  =-\Omega(X_{\varphi},X_{\psi})\\
&  =\Omega(X_{\psi},X_{\varphi})\\
&  =(i_{X_{\psi}}\Omega)(X_{\varphi})\\
&  =\left[  \widetilde{d}(\psi)\right]  (X_{\varphi})\\
&  =\widetilde{X_{\varphi}}(\psi)
\end{align*}
et
\begin{align*}
i_{[X_{\varphi},X_{\psi}]}\Omega &  =[\theta_{X_{\varphi}},i_{X_{\psi}%
}](\Omega)\\
&  =\theta_{X_{\varphi}}(i_{X_{\psi}}\Omega)-i_{X_{\psi}}(\theta_{X_{\varphi}%
}\Omega)\\
&  =\theta_{X_{\varphi}}\widetilde{d}\psi\\
&  =\widetilde{d}[\theta_{X_{\varphi}}(\psi)]\\
&  =\widetilde{d}[\widetilde{X_{\varphi}}(\psi)]\\
&  =\widetilde{d}\{\varphi,\psi\}\\
&  =i_{\{\varphi,\psi\}}\Omega\text{.}%
\end{align*}
Comme $\Omega$ est nond\'{e}g\'{e}n\'{e}r\'{e}e, on d\'{e}duit que
\[
\lbrack X_{\varphi},X_{\psi}]=X_{\{\varphi,\psi\}}\text{.}%
\]

D'o\`{u} les deux assertions.
\end{demonstration}

\begin{theoreme}
Si le couple $(\mathcal{D}(M^{A}),\mathbf{\sim})$ admet une structure de
$A$-alg\`{e}bre de Lie-Rinehart-Jacobi symplectique, alors l'application
\[
C^{\infty}(M^{A},A)\times C^{\infty}(M^{A},A)\longrightarrow C^{\infty}%
(M^{A},A),(\varphi,\psi)\longmapsto\left\{  \varphi,\psi\right\}  ,
\]
d\'{e}finit une structure de $A$-alg\`{e}bre de Jacobi sur $C^{\infty}%
(M^{A},A)$.
\end{theoreme}

\begin{demonstration}
\bigskip\ Comme, pour tous $\varphi,\psi\in C^{\infty}(M^{A},A)$, on a%
\begin{align*}
\left\{  \varphi,\psi\right\}   &  =-\Omega(X_{\varphi},X_{\psi})\\
&  =\Omega(X_{\psi},X_{\varphi})\\
&  =-\left[  -\Omega(X_{\psi},X_{\varphi})\right] \\
&  =-\{\psi,\varphi\}\text{.}%
\end{align*}
D'o\`{u}%
\[
\{\varphi,\ \psi\}=-\{\psi,\varphi\}\text{.}%
\]
$\ $Pour tout $a\in A$,
\begin{align*}
\{\varphi,a\cdot\psi\}  &  =\widetilde{X_{\varphi}}(a\psi)\\
&  =a\cdot\{\varphi,\ \psi\}\text{.}%
\end{align*}

Ainsi
\[
\{\varphi,a\cdot\psi\}=a\cdot\{\varphi,\ \psi\}\text{.}%
\]

Pour $\varphi,\psi,\nu\in C^{\infty}(M^{A},A)$, on a%
\begin{align*}
&  \{\varphi,\{\psi,\nu\}\}+\{\psi,\{\ \nu,\varphi\}\}+\{\nu,\{\varphi
,\psi\}\}\\
&  =\{\varphi,\{\psi,\nu\}\}-\{\psi,\{\ \varphi,\nu\}\}-\{\{\varphi,\psi
\},\nu\}\\
&  =\widetilde{X_{\varphi}}[\widetilde{X_{\psi}}(\nu)]-\widetilde{X_{\psi}%
}[\widetilde{X_{\varphi}}(\nu)]-\ \widetilde{X_{\{\varphi,\psi\}}}(\nu)\\
&  =[\widetilde{X_{\varphi}},\widetilde{X_{\psi}}](\nu)-\widetilde
{[X_{\varphi},X_{\psi}]}(\nu)\\
&  =\left(  [\widetilde{X_{\varphi}},\widetilde{X_{\psi}}]\ -\widetilde
{[X_{\varphi},X_{\psi}]}\right)  (\nu)\\
&  =0.
\end{align*}
L'identit\'{e} de Jacobi est ainsi d\'{e}montr\'{e}e.

Pour $\varphi\in C^{\infty}(M^{A},A)$, l'application
\[
ad(\varphi):C^{\infty}(M^{A},A)\longrightarrow C^{\infty}(M^{A},A),\psi
\longmapsto\{\varphi,\psi\},
\]
est un op\'{e}rateur diff\'{e}rentiel. En effet comme
\[
\widetilde{X_{\varphi}}:C^{\infty}(M^{A},A)\longrightarrow C^{\infty}%
(M^{A},A)
\]
est un op\'{e}rateur diff\'{e}rentiel, alors pour $\psi_{1},\psi_{2}\in
C^{\infty}(M^{A},A)$ on a:%
\begin{align*}
&  \left[  ad(\varphi)\right]  (\psi_{1}\cdot\psi_{2})\\
&  =\{\varphi,\psi_{1}\cdot\psi_{2}\}\\
&  =\widetilde{X_{\varphi}}(\psi_{1}\cdot\psi_{2})\\
&  =\widetilde{X_{\varphi}}(\psi_{1})\cdot\psi_{2}+\psi_{1}\cdot
\widetilde{X_{\varphi}}(\psi_{2})-\psi_{1}\cdot\psi_{2}\cdot\widetilde
{X_{\varphi}}(1_{C^{\infty}(M^{A},A)})\\
&  =\{\varphi,\psi_{1}\}\cdot\psi_{2}+\psi_{1}\cdot\{\varphi,\psi_{2}%
\}-\psi_{1}\cdot\psi_{2}\cdot\{\varphi,1_{C^{\infty}(M^{A},A)}\}\\
&  =\left[  ad(\varphi)\right]  (\psi_{1})\cdot\psi_{2}+\psi_{1}\cdot\left[
ad(\varphi)\right]  (\psi_{2})-\psi_{1}\cdot\psi_{2}\cdot\left[
ad(\varphi)\right]  (1_{C^{\infty}(M^{A},A)})\text{.}%
\end{align*}
On conclut que $C^{\infty}(M^{A},A)$ est une $A$-alg\`{e}bre de Jacobi
c'est-\`{a}-dire que $M^{A}$ est une $A$-vari\'{e}t\'{e} de Jacobi.
\end{demonstration}

\subsection{$A$-Structure de Jacobi sur $M^{A}$ lorsque $M$ est une
vari\'{e}t\'{e} de Jacobi}

On rappelle qu'une structure de vari\'{e}t\'{e} de Jacobi sur une
vari\'{e}t\'{e} lisse $M$ est la donn\'{e}e d'une structure d'alg\`{e}bre de
Lie r\'{e}elle sur $C^{\infty}(M)$, de crochet, $\left\{  ,\right\}  $, telle
que pour tout $f\in C^{\infty}(M)$, l'application
\[
ad(f):C^{\infty}(M)\longrightarrow C^{\infty}(M),g\mapsto\left\{  f,g\right\}
,
\]
$\ $soit un op\'{e}rateur diff\'{e}rentiel. Dans ce cas, on dit que $M$ est
une vari\'{e}t\'{e} de Jacobi et que $C^{\infty}(M)$ est une alg\`{e}bre de Jacobi.

Dans ces conditions, pour tous $f,g\in C^{\infty}(M)$, on a
\[
ad(fg)=f\cdot ad(g)+g\cdot ad(f)-f\cdot g\cdot ad(1_{C^{\infty}(M)})\text{.}%
\]
Lorsque $ad(1_{C^{\infty}(M)})=0$, on dit que $M$ est une vari\'{e}t\'{e} de Poisson.

Pour tout $f\in C^{\infty}(M)$,
\[
\left[  ad(f)\right]  ^{A}:C^{\infty}(M)\longrightarrow C^{\infty}%
(M^{A},A),g\longmapsto\left\{  f,g\right\}  ^{A},
\]
est un op\'{e}rateur diff\'{e}rentiel et
\[
\widetilde{\left[  ad(f)\right]  ^{A}}:C^{\infty}(M^{A},A)\longrightarrow
C^{\infty}(M^{A},A)
\]
est l'unique op\'{e}rateur diff\'{e}rentiel qui est $A$-lin\'{e}aire et qui
est tel que
\begin{align*}
\widetilde{\left[  ad(f)\right]  ^{A}}(g^{A})  &  =\left[  ad(f)\right]
^{A}(g)\\
&  =\left\{  f,g\right\}  ^{A}%
\end{align*}
pour tout $g\in C^{\infty}(M)$.

\begin{proposition}
Pour $\varphi\in C^{\infty}(M^{A},A)$, l'application
\[
\tau_{\varphi}:C^{\infty}(M)\longrightarrow C^{\infty}(M^{A},A),f\longmapsto
-\widetilde{[ad(f)]^{A}}(\varphi),
\]
est un op\'{e}rateur diff\'{e}rentiel.
\end{proposition}

\begin{demonstration}
L'application $\tau_{\varphi}$ est manifestement lin\'{e}aire. Pour $f,g\in
C^{\infty}(M)$, on a%
\begin{align*}
\tau_{\varphi}(fg)  &  =-\widetilde{[ad(fg)]^{A}}(\varphi)\\
&  =-\widetilde{[f\cdot ad(g)+g\cdot ad(f)-f\cdot g\cdot ad(1_{C^{\infty}%
(M)})]^{A}}(\varphi)\\
&  =f^{A}\cdot\left(  -\widetilde{[ad(g)]^{A}}\right)  (\varphi)+g^{A}%
\cdot\left(  -\widetilde{[ad(f)]^{A}}\right)  (\varphi)-f^{A}\cdot g^{A}%
\cdot\left(  -\widetilde{\left[  ad(1_{C^{\infty}(M)})\right]  ^{A}}\right)
(\varphi)\\
&  =\left(  -\widetilde{[ad(f)]^{A}}\right)  (\varphi)\cdot g^{A}+f^{A}%
\cdot\left(  -\widetilde{[ad(g)]^{A}}\right)  (\varphi)-f^{A}\cdot g^{A}%
\cdot\left(  -\widetilde{\left[  ad(1_{C^{\infty}(M)})\right]  ^{A}}\right)
(\varphi)\\
&  =\tau_{\varphi}(f)\cdot g^{A}+f^{A}\cdot\tau_{\varphi}(g)-f^{A}\cdot
g^{A}\cdot\tau_{\varphi}(1_{C^{\infty}(M)})\text{.}%
\end{align*}
D'o\`{u} l'assertion.
\end{demonstration}

Pour $\varphi\in C^{\infty}(M^{A},A)$, l'application
\[
\widetilde{\tau_{\varphi}}:C^{\infty}(M^{A},A)\longrightarrow C^{\infty}%
(M^{A},A)
\]
est l'unique op\'{e}rateur diff\'{e}rentiel qui est $A$-lin\'{e}aire et qui
est tel que
\[
\widetilde{\tau_{\varphi}}(f^{A})=\tau_{\varphi}(f)
\]
pour tout $f\in C^{\infty}(M)$.

\begin{theoreme}
Si $M$ est une vari\'{e}t\'{e} de Jacobi, de crochet $\left\{  ,\right\}  $,
alors l'application
\[
\left\{  ,\right\}  _{A}:C^{\infty}(M^{A},A)\times C^{\infty}(M^{A}%
,A)\longrightarrow C^{\infty}(M^{A},A),(\varphi,\psi)\longmapsto
\widetilde{\tau}_{\varphi}(\psi),\text{ }%
\]
d\'{e}finit une $A$-structure de Jacobi sur $M^{A}$.
\end{theoreme}

\begin{demonstration}
Les techniques de d\'{e}monstration sont identiques \`{a} celles utilis\'{e}es
dans \cite{bo2}.
\end{demonstration}

On dit que la $A$-structure de Jacobi sur $M^{A}$\ d\'{e}finie par $\left\{
,\right\}  _{A}$ est le prolongement \`{a} $M^{A}$\ de la structure de
Jacobi\ sur $M$ d\'{e}finie par $\left\{  ,\right\}  $.

\subsection{$A$-Structure de Jacobi sur $M^{A}$ lorsque $M$ est une
vari\'{e}t\'{e} localement conform\'{e}ment symplectique}

\begin{proposition}
\bigskip Si
\[
\omega:\mathfrak{X}(M)\times\mathfrak{X}(M)\longrightarrow C^{\infty}(M)
\]
est une $2$-forme diff\'{e}rentielle nond\'{e}g\'{e}n\'{e}r\'{e}e, alors
\[
\omega^{A}:\mathfrak{X}(M^{A})\times\mathfrak{X}(M^{A})\longrightarrow
C^{\infty}(M^{A},A)
\]
est nond\'{e}g\'{e}n\'{e}r\'{e}e.
\end{proposition}

\begin{demonstration}
On va v\'{e}rifier que l'application
\[
\mathfrak{X}(M^{A})\longrightarrow\wedge^{1}\left(  M^{A},A\right)
,X\longmapsto i_{X}\omega^{A},
\]
est un isomorphisme de $C^{\infty}(M^{A},A)$-modules. Pour cela, il suffit de
montrer qu'en chaque point $\xi\in M^{A}$, l'application%
\[
T_{\xi}M^{A}\longrightarrow T_{\xi}^{\ast}M^{A},v\longmapsto i_{v}\omega
^{A}\left(  \xi\right)  ,
\]
est un isomorphisme de $A$-modules.

On note $\mathfrak{m}$ l'unique id\'{e}al maximal de $A$ et $h$ la hauteur de
$A$: $h$ est l'entier naturel tel que $\mathfrak{m}^{h}\neq(0)$ et
$\mathfrak{m}^{h+1}=(0)$. Il existe des sous-espaces vectoriels $V_{0}%
,V_{1},...,V_{h}$ de $A$ tels que
\[
A=V_{0}\oplus...\oplus V_{h}%
\]
avec $V_{i}\subset\mathfrak{m}^{i}$ et $V_{i}\oplus\mathfrak{m}^{i+1}%
=\mathfrak{m}^{i}$ pour $i=0,1,...,h$. Ainsi, $V_{0}=%
\mathbb{R}
$ et $V_{h}=\mathfrak{m}^{h}$.

Pour $k=0,1,...,h$, on note $I_{k}$ l'ensemble des indices d'une base de
$V_{k}$ et $\left(  a_{\alpha_{k}}\right)  _{_{\alpha_{k}}\in I_{k}}$ une base
de $V_{k}$. Soit $\xi\in M^{A}$, $x_{0}\in M$ l'origine de $\xi$ et $\left(
U,\varphi\right)  $ une carte locale de $M$ en $x_{0}$, de fonctions
coordonn\'{e}es $\left(  x_{1},...,x_{n}\right)  $ o\`{u} $n=\dim M$.

Injection : \ Soit $v\in T_{\xi}M^{A}$ tel que $\omega^{A}\left(  \xi\right)
(v,w)=0$ pour tout $w\in T_{\xi}M^{A}$. L'espace tangent $T_{\xi}M^{A}$ est un
$A$-module libre de rang $n$ dont une base est $\left(  \left(  \frac
{\partial}{\partial x_{1}}\right)  ^{A}|_{\xi},...,\left(  \frac{\partial
}{\partial x_{n}}\right)  ^{A}|_{\xi}\right)  $. En particulier, on a%
\[
\omega^{A}\left(  \xi\right)  (v,\left(  \frac{\partial}{\partial x_{l}%
}\right)  ^{A}|_{\xi})=0
\]
pour $l=1,2,...,n$.

On a
\[
v=\underset{j=1}{\overset{n}{\sum}}v_{j}\cdot\left(  \frac{\partial}{\partial
x_{j}}\right)  ^{A}|_{\xi}%
\]
avec $v_{j}\in A$, pour $j=1,2,...,n$. Comme%
\[
v_{j}=\sum_{\alpha_{k}\in I_{k},k=0,1,...,h}v_{j,\alpha_{k}}\cdot
a_{\alpha_{k}}%
\]
avec $v_{j,\alpha_{k}}\in%
\mathbb{R}
$, alors%
\[
v=\underset{j=1}{\overset{n}{\sum}}\left[  \sum_{\alpha_{k}\in I_{k}%
,k=0,1,...,h}v_{j,\alpha_{k}}\cdot a_{\alpha_{k}}\right]  \cdot\left(
\frac{\partial}{\partial x_{j}}\right)  ^{A}|_{\xi}\text{.}%
\]
L'\'{e}quation%
\[
\omega^{A}\left(  \xi\right)  (v,\left(  \frac{\partial}{\partial x_{l}%
}\right)  ^{A}|_{\xi})=0
\]
signifie que
\[
\sum_{\alpha_{k}\in I_{k},k=0,1,...,h}\left(  \sum\limits_{j=1}^{n}%
v_{j,\alpha_{k}}\cdot\left[  \omega\left(  \frac{\partial}{\partial x_{j}%
},\frac{\partial}{\partial x_{l}}\right)  \right]  ^{A}(\xi)\right)  \cdot
a_{\alpha_{k}}=0
\]
pour $l=1,2,...,n$. Ainsi, on a
\begin{align*}
&  \sum_{\alpha_{k}\in I_{k},k=0,1,...,h}\left(  \sum\limits_{j=1}%
^{n}v_{j,\alpha_{k}}\cdot\left[  \omega\left(  \frac{\partial}{\partial x_{j}%
},\frac{\partial}{\partial x_{l}}\right)  \right]  (x_{0})\right)  \cdot
a_{\alpha_{k}}\\
&  +\sum_{\alpha_{k}\in I_{k},k=0,1,...,h}\left(  \sum\limits_{j=1}%
^{n}v_{j,\alpha_{k}}\cdot\theta_{jl}(\xi)\right)  \cdot a_{\alpha_{k}}\\
&  =0
\end{align*}
avec $\theta_{jl}(\xi)\in\mathfrak{m}$. Il s'ensuit que $\sum\limits_{j=1}%
^{n}v_{j,\alpha_{0}}\cdot\left[  \omega\left(  \frac{\partial}{\partial x_{j}%
},\frac{\partial}{\partial x_{l}}\right)  \right]  (x_{0})=0$ pour
$l=1,2,...,n$. Comme $\omega$ est nond\'{e}g\'{e}n\'{e}r\'{e}e, on d\'{e}duit
que $v_{j,\alpha_{0}}=0$ pour $j=1,2,...,n$.

En raisonnant par r\'{e}currence, on suppose que $v_{j,\alpha_{r}}=0$ pour
$j=1,2,...,n$. Montrons que $v_{j,\alpha_{r+1}}=0$ pour $j=1,2,...,n$.
L'\'{e}quation devient%
\[
\sum_{\alpha_{k}\in I_{k},k=r+1,...,h}\left(  \sum\limits_{j=1}^{n}%
v_{j,\alpha_{k}}\cdot\left[  \omega\left(  \frac{\partial}{\partial x_{j}%
},\frac{\partial}{\partial x_{l}}\right)  \right]  ^{A}(\xi)\right)  \cdot
a_{\alpha_{k}}=0
\]

pour $l=1,2,...,n$. Ainsi
\begin{align*}
0  &  =\sum_{\alpha_{r+1}\in I_{r+1}}\left(  \sum\limits_{j=1}^{n}%
v_{j,\alpha_{r+1}}\cdot\left[  \omega\left(  \frac{\partial}{\partial x_{j}%
},\frac{\partial}{\partial x_{l}}\right)  \right]  ^{A}(\xi)\right)  \cdot
a_{\alpha_{r+1}}\\
&  +\sum_{\alpha_{k}\in I_{k},k=r+2,...,h}\left(  \sum\limits_{j=1}%
^{n}v_{j,\alpha_{k}}\cdot\left[  \omega\left(  \frac{\partial}{\partial x_{j}%
},\frac{\partial}{\partial x_{l}}\right)  \right]  ^{A}(\xi)\right)  \cdot
a_{\alpha_{k}}\text{.}%
\end{align*}
On a%
\begin{align*}
0  &  =\sum_{\alpha_{r+1}\in I_{r+1}}\left(  \sum\limits_{j=1}^{n}%
v_{j,\alpha_{r+1}}\cdot\left[  \omega\left(  \frac{\partial}{\partial x_{j}%
},\frac{\partial}{\partial x_{l}}\right)  \right]  (x_{0})+\theta_{ij}%
(\xi)\right)  \cdot a_{\alpha_{r+1}}\\
&  +\sum_{\alpha_{k}\in I_{k},k=r+2,...,h}\left(  \sum\limits_{j=1}%
^{n}v_{j,\alpha_{k}}\cdot\left[  \omega\left(  \frac{\partial}{\partial x_{j}%
},\frac{\partial}{\partial x_{l}}\right)  \right]  ^{A}(\xi)\right)  \cdot
a_{\alpha_{k}}\text{,}%
\end{align*}
avec $\theta_{ij}(\xi)\in\mathfrak{m}$. D'o\`{u}
\[
\sum_{\alpha_{r+1}\in I_{r+1}}\left(  \sum\limits_{j=1}^{n}v_{j,\alpha_{r+1}%
}\cdot\left[  \omega\left(  \frac{\partial}{\partial x_{j}},\frac{\partial
}{\partial x_{l}}\right)  \right]  (x_{0})\right)  \cdot a_{\alpha_{r+1}%
}=0\text{.}%
\]
Comme $a_{\alpha_{r+1}}$ est une base de $V_{r+1}$, on d\'{e}duit
\[
\sum\limits_{j=1}^{n}v_{j,\alpha_{r+1}}\cdot\left[  \omega\left(
\frac{\partial}{\partial x_{j}},\frac{\partial}{\partial x_{l}}\right)
\right]  (x_{0})=0
\]
pour $l=1,2,...,n$. Comme $\omega$ est nond\'{e}g\'{e}n\'{e}r\'{e}e, alors
$v_{j,\alpha_{r+1}}=0$.

On conclut que $v=0$.

Surjection : L'espace $T_{\xi}^{\ast}M^{A}$ est un $A$-module libre de rang
$n=\dim M$. Une base est $\left(  \left(  dx_{1}\right)  ^{A}|_{\xi
},...,\left(  dx_{n}\right)  ^{A}|_{\xi}\right)  $ o\`{u} $\left(
x_{1},...,x_{n}\right)  $ est un syst\`{e}me de coordonn\'{e}es locales au
voisinage $U$ de l'origine $x_{0}$ de $\xi\in M^{A}$.

Soit $\eta\in$ $T_{\xi}^{\ast}M^{A}$. On \'{e}crit
\[
\eta=\underset{k=1}{\overset{n}{\sum}}\eta_{k}\cdot\left(  dx_{k}\right)
^{A}|_{\xi}%
\]
avec $\eta_{k}\in A$ pour $k=1,2,...,n$.

Soit $\theta_{k}$ l'unique champ de vecteurs sur $U$ tel que
\[
dx_{k}=i_{\theta_{k}}\omega|_{U}\text{.}%
\]
Alors
\[
\left(  dx_{k}\right)  ^{A}|_{\xi}=i_{\theta_{k}^{A}(\xi)}\cdot\omega^{A}%
(\xi)\text{.}%
\]
Ainsi,
\[
\eta=i_{v}\omega^{A}(\xi)
\]
avec
\[
v=\underset{i=1}{\overset{n}{\sum}}\eta_{i}\cdot\theta_{i}^{A}(\xi)\in T_{\xi
}M^{A}\text{.}%
\]

Ce qui ach\`{e}ve la d\'{e}monstration.
\end{demonstration}

Dans toute la suite $(M,\alpha,\omega)$ d\'{e}signe une vari\'{e}t\'{e}
localement conform\'{e}ment symplectique de $1$-forme $\alpha$ et de $2$-forme
$\omega$. Dans ce cas, l'application%
\[
\rho_{\alpha}(\theta):C^{\infty}(M)\longrightarrow C^{\infty}(M),f\longmapsto
\theta(f)+f\cdot\alpha(\theta),
\]
est un op\'{e}rateur diff\'{e}rentiel pour tout $\theta\in\mathfrak{X}(M)$.

De plus l'application%

\[
\rho_{\alpha}:\mathfrak{X}(M)\longrightarrow\mathcal{D}(M),\theta
\longmapsto\rho_{\alpha}(\theta),
\]
est une repr\'{e}sentation et le triplet $(\mathfrak{X}(M),\rho_{\alpha
},\omega)$ est une alg\`{e}bre de Lie-Rinehart-Jacobi \cite{ok3},\cite{ok4} .
L'op\'{e}rateur de cohomologie associ\'{e}e \`{a} la repr\'{e}sentation
$\rho_{\alpha}$ est l'op\'{e}rateur de cohomologie de Lichnerowicz $d_{\alpha
}$ \cite{lic}.

La vari\'{e}t\'{e} diff\'{e}rentielle $M$ est une vari\'{e}t\'{e} de Jacobi
o\`{u} le crochet de deux fonctions $f$ et $g$ est donn\'{e} par
\[
\left\{  f,g\right\}  =-\omega(X_{f},X_{g}),
\]
$X_{f}$ \'{e}tant l'unique champ de vecteurs sur $M$ tel que
\[
i_{X_{f}}\omega=d_{\alpha}f\text{.}%
\]

Pour $X\in\mathfrak{X}(M^{A})$, consid\'{e}r\'{e} comme d\'{e}rivation de
$C^{\infty}(M)$ dans $C^{\infty}(M^{A},A)$, l'application%
\[
\rho_{\alpha^{A}}(X):C^{\infty}(M^{A},A)\longrightarrow C^{\infty}%
(M^{A},A),\varphi\longmapsto\widetilde{X}(\varphi)+\varphi\cdot\alpha^{A}(X),
\]
est un op\'{e}rateur diff\'{e}rentiel qui est $A$-lin\'{e}aire et
l'application
\[
\rho_{\alpha^{A}}:\mathfrak{X}(M^{A})\longrightarrow\mathcal{D}_{A}\left[
C^{\infty}(M^{A},A)\right]  ,X\longmapsto\rho_{\alpha^{A}}(X),
\]
est $C^{\infty}(M^{A},A)$-lin\'{e}aire et est un morphisme de $A$-alg\`{e}bres
de Lie.

\begin{proposition}
Le couple $(\mathfrak{X}(M^{A}),\rho_{\alpha^{A}})$ est une $A$-alg\`{e}bre de Lie-Rinehart.
\end{proposition}

La d\'{e}monstration est une simple v\'{e}rification.

Si $d_{\rho_{\alpha^{A}}}$ d\'{e}signe l'op\'{e}rateur de cohomologie
associ\'{e} \`{a} la repr\'{e}sentation $\rho_{\alpha^{A}}$ et si $d^{A}$,
\cite{bos}, est l'op\'{e}rateur de cohomologie associ\'{e} \`{a} la
repr\'{e}sentation
\[
\mathfrak{X}(M^{A})\longrightarrow Der_{A}\left[  C^{\infty}(M^{A},A)\right]
,X\longmapsto\widetilde{X,}%
\]
on v\'{e}rifie que
\[
d_{\rho_{\alpha^{A}}}\eta=d^{A}\eta+\alpha^{A}\Lambda\eta
\]
pour tout $\eta\in\Lambda(M^{A},A)$. Ainsi on conclut que
\[
d_{\rho_{\alpha^{A}}}=d_{\alpha^{A}}^{A}\text{.}%
\]
On a imm\'{e}diatement:

\begin{proposition}
Si $\left(  M,\alpha,\omega\right)  $ est une vari\'{e}t\'{e} localement
conform\'{e}ment symplectique, alors le triplet $\left(  \mathfrak{X}%
(M^{A}),\rho_{\alpha^{A}},\omega^{A}\right)  $ est une $A$-alg\`{e}bre de
Lie-Rinehart-Jacobi symplectique.
\end{proposition}

\begin{demonstration}
On a
\begin{align*}
d_{\rho_{\alpha^{A}}}\omega^{A}  &  =d_{\alpha^{A}}^{A}\omega^{A}\\
&  =d^{A}\omega^{A}+\alpha^{A}\Lambda\omega^{A}\\
&  =(d\omega+\alpha\Lambda\omega)^{A}\\
&  =0\text{.}%
\end{align*}
Comme le couple $\mathfrak{X}(M^{A}),\rho_{\alpha^{A}})$ est une
$A$-alg\`{e}bre de Lie-Rinehart, comme $\omega^{A}$ est
nond\'{e}g\'{e}n\'{e}r\'{e}e et comme $d_{\alpha^{A}}^{A}\omega^{A}=0$, on
conclut que le triplet $\left(  \mathfrak{X}(M^{A}),\rho_{\alpha^{A}}%
,\omega^{A}\right)  $ est une $A$-alg\`{e}bre de Lie-Rinehart-Jacobi symplectique.
\end{demonstration}

Pour $F\in C^{\infty}(M^{A},A)$, on note $X_{F}$ l'unique \'{e}l\'{e}ment de
$\mathfrak{X}(M^{A})$ tel que
\[
i_{X_{F}}\omega^{A}=d_{\alpha^{A}}^{A}F\text{.}%
\]

En s'inpirant des techniques de \cite{ok3}, page $1087$, on d\'{e}duit:

\begin{proposition}
L'application
\[
\left\{  ,\right\}  _{\omega^{A}}:C^{\infty}(M^{A},A)\times C^{\infty}%
(M^{A},A)\longrightarrow C^{\infty}(M^{A},A),(F,G)\longmapsto-\omega^{A}%
(X_{F},X_{G}),
\]
d\'{e}finit une structure de $A$-alg\`{e}bre de Jacobi sur $C^{\infty}%
(M^{A},A)$.
\end{proposition}

\begin{proposition}
Pour $f\in C^{\infty}(M)$,
\[
X_{f^{A}}=(X_{f})^{A}\text{.}%
\]

\end{proposition}

\begin{demonstration}
On a
\begin{align*}
i_{X_{f^{A}}}\omega^{A}  &  =d_{\alpha^{A}}^{A}f^{A}\\
&  =d^{A}f^{A}+f^{A}\cdot\alpha^{A}\\
&  =(d_{\alpha}f)^{A}\\
&  =(i_{X_{f}}\omega)^{A}\\
&  =i_{(X_{f})^{A}}\omega^{A}\text{.}%
\end{align*}
Comme $\omega^{A}$ est nond\'{e}g\'{e}n\'{e}r\'{e}e, l'assertion s'ensuit.
\end{demonstration}

\begin{proposition}
Si $\left\{  ,\right\}  $ est le crochet de Jacobi d\'{e}fini sur $C^{\infty
}(M)$ par la structure de vari\'{e}t\'{e} de Jacobi d\'{e}duite de la
vari\'{e}t\'{e} localement conform\'{e}ment symplectique $\left(
M,\alpha,\omega\right)  $, alors pour $f,g\in C^{\infty}(M)$ on a
\[
\left\{  f^{A},g^{A}\right\}  _{\omega^{A}}=\left\{  f,g\right\}  ^{A}\text{.}%
\]

\end{proposition}

\begin{demonstration}
On a
\begin{align*}
\left\{  f^{A},g^{A}\right\}  _{\omega^{A}}  &  =-\omega^{A}(X_{f^{A}%
},X_{g^{A}})\\
&  =-\omega^{A}((X_{f})^{A},(X_{g})^{A})\\
&  =\left[  -\omega(X_{f},X_{g})\right]  ^{A}\\
&  =\left\{  f,g\right\}  ^{A}\text{.}%
\end{align*}
D'o\`{u} l'assertion.
\end{demonstration}

La proposition pr\'{e}c\'{e}dente signifie que si $\left(  M,\alpha
,\omega\right)  $ est une vari\'{e}t\'{e} localement conform\'{e}ment
symplectique, la $A$-structure d'alg\`{e}bre de Jacobi sur $C^{\infty}%
(M^{A},A)$ d\'{e}finie par la $A$-alg\`{e}bre de Lie-Rinehart-Jacobi
symplectique $\left(  \mathfrak{X}(M^{A}),\rho_{\alpha^{A}},\omega^{A}\right)
$ co\"{\i}ncide avec le prolongement \`{a} $M^{A}$ de la structure de Jacobi
sur $M$ d\'{e}finie par la vari\'{e}t\'{e} localement conform\'{e}ment
symplectique $\left(  M,\alpha,\omega\right)  $.

\end{document}